\documentclass[11pt]{article}
\PassOptionsToPackage{obeyspaces}{url}
\usepackage[colorlinks=true,citecolor=blue,urlcolor=blue,linkcolor=blue,bookmarksopen=true]{hyperref}
\usepackage{amsfonts, amsmath, amssymb, amsthm}

\usepackage{breakurl}

\usepackage{etoolbox} 
\patchcmd{\thebibliography}{\leftmargin\labelwidth}{\leftmargin\labelwidth\addtolength\itemsep{-0.2\baselineskip}}{}{}

\usepackage{mathtools}

\author{Boris Bukh\thanks{Department of Mathematical Sciences, Carnegie Mellon University, Pittsburgh, PA 15213, USA. Supported in part by Sloan Research Fellowship and by U.S.\ taxpayers through NSF CAREER grant DMS-1555149.} \and Oleksandr Rudenko\thanks{Department of Mathematical Sciences, Carnegie Mellon University, Pittsburgh, PA 15213, USA. Supported in part by U.S.\ taxpayers through NSF CAREER grant DMS-1555149.}}

\title{Order-isomorphic twins in permutations}
\date{}

\usepackage[nameinlink]{cleveref}
\newtheorem{theorem}{Theorem}
\newtheorem{lemma}[theorem]{Lemma}

\newtheorem{claim}[theorem]{Claim}

\newcommand*{\eqdef}{\stackrel{\mbox{\normalfont\tiny def}}{=}}   
\DeclarePairedDelimiter\abs{\lvert}{\rvert}                     

\newcommand*{\perm}{\mathbf}                                    

\begin{document}
\maketitle

\begin{abstract}
Let $a_1,\dotsc,a_n$ be a permutation of $[n]$. Two disjoint order-isomorphic subsequences are called \emph{twins}.
We show that every permutation of $[n]$ contains twins of length $\Omega(n^{3/5})$ improving the trivial bound
of $\Omega(n^{1/2})$. We also show that a random permutation contains twins of length $\Omega(n^{2/3})$, which is sharp.
\end{abstract}

In this paper we regard permutations as sequences of symbols, devoid of any group-theoretic meaning.
So, for us a \emph{permutation} on a finite set $\Sigma$ is a sequence of elements of $\Sigma$ in which each element of $\Sigma$
appears exactly once. We call a subsequence of a permutation \emph{subpermutation}. For instance, $135642$ is a permutation
of $[6]$, and $1562$ is a subpermutation inside, which itself is a permutation of $\{1,2,5,6\}$. We denote permutations
by bold letters.

Throughout the paper we consider only the permutations of finite sets of natural numbers.
We say that permutations $\perm{a}=(a_1,\dotsc,a_L)$ and $\perm{b}=(b_1,\dotsc,b_L)$ are \emph{order-isomorphic}
if $(a_i<a_j) \iff (b_i<b_j)$. For example, $1562$ is order-isomorphic to $1342$.

We call a pair of subpermutation $\perm{a},\perm{b}$ of $\perm{c}$ \emph{twins} if $\perm{a}$ and $\perm{b}$ are order-isomorphic
and disjoint (do not contain the same symbol). For example, $152$ and $364$ are twins in $135642$ of length~$3$.
We denote by $t(n)$ the largest integer such that every permutation of $[n]$ contains a pair of twins of length~$t(n)$.

The problem of estimating $t(n)$ was raised by Gawron \cite{gawron}, who observed that $t(n)\geq (n^{1/2}-1)/2$ follows from the Erd\H{o}s--Szekeres
theorem, and that $t(n)=O(n^{2/3})$ follows from the first moment method. He further conjectured that $t(n)=\Omega(n^{2/3})$. This is not known even
for random permutations: the best result is due to Dudek, Grytczuk, and Ruci\'nski \cite{dgr} who showed that
a random permutation almost surely contains twins of length $\Omega(n^{2/3}/\log^{1/3} n)$.

In this short note, we give a first non-trivial lower bound on $t(n)$, and remove the logarithmic factor from the Dudek--Grytczuk--Ruci\'nski result.
\begin{theorem}\label{thm:main}
For $n\geq 2$, every permutation of $[n]$ contains twins of length at least $\tfrac{1}{8}n^{3/5}$.
\end{theorem}
\begin{theorem}\label{thm:random}
A random permutation of $[n]$ almost surely contains twins of length at least $\tfrac{1}{80}n^{2/3}$, as $n\to\infty$.
\end{theorem}
In view of Gawron's result, \Cref{thm:random} is sharp up to the constant factor.

\section*{Proof of \Cref{thm:main}}
The proof relies on 
a result of Beame and Huynh-Ngoc \cite[Lemma 5.9]{bhn}, which previously 
was used by Bukh and Zhou \cite{bukh_zhou} to study a related notion of twins in words.
\begin{lemma}\label{lem:bhn}
Among any three permutations $\perm{c}^{(0)},\perm{c}^{(1)},\perm{c}^{(2)}$ of $[m]$ we may find two distinct, say, $\perm{c}^{(k)}$ and $\perm{c}^{(\ell)}$,
that contain the same subpermutation of length at least $m^{1/3}$.
\end{lemma}

Call twin subpermutations $a_1,\dotsc,a_L$ and $b_1,\dotsc,b_L$ \emph{close} if $\abs{b_i-a_i}\leq n^{2/5}$ for all $i$.
Let $t'(m,n)$ be the largest integer so that whenever $\Sigma\subset [n]$ is any set of at least $m$ elements, every
permutation of $\Sigma$ contains close twins of length at least $t'(m,n)$. 
\begin{claim}\label{claim}
If $m\geq 7n^{3/5}$, then $t'(m,n)\geq t'(m-7n^{3/5},n)+n^{1/5}$.
\end{claim}
From $\lfloor \tfrac{1}{7}n^{2/5}\rfloor$ many invocations of \Cref{claim} we infer that $t(n)\geq t'(n,n)\geq \lfloor \tfrac{1}{7}n^{2/5}\rfloor n^{1/5}$,
implying \Cref{thm:main} for $n\geq 56^{5/2}$. When $2\leq n\leq 56^{5/2}$, \Cref{thm:main} follows from
$t(n)\geq (n^{1/2}-1)/2\geq \tfrac{1}{8}n^{3/5}$.\medskip

We now prove the claim. We can clearly assume that $\abs{\Sigma}=m$.
Let $\perm{a}=(a_1,\dots,a_m)$ be an arbitrary permutation of $\Sigma$. Consider its first $3r$ elements, where $r\eqdef \lceil 2n^{3/5}\rceil$. 
Say $a_1,\dotsc,a_{3r}$ is a permutation of the set $\{b_0,\dotsc,b_{3r-1}\}$, where
$b_0<\dotsb<b_{3r-1}$. Consider the triples $(b_0,b_1,b_2)$, $(b_3,b_4,b_5)$, $\dots$, $(b_{3r-3},b_{3r-2},b_{3r-1})$.
Since $\sum_{i=0}^{r} (b_{3i+2}-b_{3i})\leq n$, the set $I_0\eqdef \{i : b_{3i+2}-b_{3i}\leq 2n/r\}$ has at least $r/2$ elements. 
For each $j=0,1,2$ let $\perm{c}^{(j)}$ be the subpermutation of $a_1,\dotsc,a_m$ obtained
by keeping only the elements $b_{3i+j}$ with $i\in I_0$. Let $c^{(j)}_i\eqdef b_{3i+j}$, and note that $\abs{c^{(j)}_i-c^{(k)}_i}\leq n^{2/5}$.

Replace each $c^{(j)}_i$ in $\perm{c}^{(j)}$ with number $i$ to obtain permutation $\tilde{\perm{c}}^{(j)}$ of~$I_0$.
By \Cref{lem:bhn} applied to the $\tilde{\perm{c}}$'s, there is $I\subset I_0$ of size $\abs{I}\geq \abs{I_0}^{1/3}\geq n^{1/5}$ and $k<\ell$ such that 
the subpermutations $\perm{c}^{(k)}_I\eqdef (c^{(k)}_i:i\in I)$ and $\perm{c}^{(\ell)}_I\eqdef(c^{(\ell)}_i:i\in I)$ are order-isomorphic. By deleting some
elements of $I$ if necessary, we may assume that $\abs{I}=\lceil n^{1/5}\rceil$. Note that
$\perm{c}^{(k)}_I$ and $\perm{c}^{(\ell)}_I$ is a pair of close twins.

Let $\perm{d}$ be the subpermutation of $\perm{a}$ obtained by deleting the first $3r$ elements, and also deleting all elements
that are contained in the intervals of the form $[\perm{c}^{(k)}_i,\perm{c}^{(k)}_i+n^{2/5}]$ for $i\in I$. Since in total
these intervals contain no more than $(n^{2/5}+1)\abs{I}$ elements, and each interval contains at least two elements among the first $3r$, 
the permutation $\perm{d}$ is of length at least $m-3r-(n^{2/5}-1)\abs{I}\geq m-7n^{3/5}$.

Given a pair of close twins $\perm{e},\perm{f}$ in $\perm{d}$, we may obtain a pair of close twins
in $\perm{a}$ by concatenating $\perm{c}^{(k)}_I$ with $\perm{e}$ and concatenating $\perm{c}^{(\ell)}_I$ with $\perm{f}$.
Indeed, let $i$ and $j$ be arbitrary, and consider two pairs of elements $\perm{c}^{(k)}_i,\perm{c}^{(\ell)}_i$ and
$e_j,f_j$. Because neither of $e_j,f_j$ is contained in the interval $T\eqdef [\perm{c}^{(k)}_i,\perm{c}^{(k)}_i+n^{2/5}]$, and
$\abs{e_j-d_j}\leq n^{2/5}$, it follows that $e_j,f_j$ are either both smaller than $\min T$ or both larger than $\max T$.
As both $\perm{c}^{(k)}_i$ and $\perm{c}^{(\ell)}_i$ are contained in $T$, we deduce that
$(\perm{c}^{(k)}_i<e_j)\iff(\perm{c}^{(\ell)}_i<f_j)$. Hence, the two concatenations indeed form a pair of twins.\vspace*{-0.4ex}

\section*{Proof of \Cref{thm:random}}
We modify the argument of Dudek--Grytczuk--Ruci\'nski. They construct a certain bipartite graph $B$ such that the matchings in $B$
correspond to twins in the original permutation. They note that $B$ contains a matching of size $v(B)/2\Delta(B)$, where $v(B)$ and $\Delta(B)$ denote the 
number of vertices and the maximum degree respectively. The logarithmic factor is lost because of the union bound to bound $\Delta(B)$.
In our proof, instead of the maximum degree, we effectively work with the typical vertex degrees. 
To help with this, we gain more independence by first Poissonizing the random process.

Let $t(\perm{p})$ be the length of the longest twin in a permutation $\perm{p}$.
We consider two ways of generating a random permutation. First, we may sample $\perm{p}$
uniformly from all permutations of $[n]$. Denote this probability distribution by $S_n$. 
Second, we may consider a Poisson process of intensity $\lambda$ on the unit square, list
the points in the order of their $x$-coordinates, and then record the relative order of $y$-coordinates.
Denote this probability distribution on permutations by $\overline{S}_{\lambda}$.

Consider an infinite sequence $p_1,p_2,\dotsc$ of independent points in $[0,1]^2$.
We may regard its prefix $p_1,\dotsc,p_m$ of length $m$ as a permutation $\perm{p}^{(m)}$ of length $m$.
We clearly have $t(\perm{p}^{(m)})\leq t(\perm{p}^{(\ell)})$ whenever $m\leq \ell$. Note that 
we may sample from $\overline{S}_n$ by sampling a number $m$ from the Poisson distribution of mean $n$
and returning $\perm{p}^{(m)}$. Since $\Pr[\operatorname{Poisson}(n/2)\geq n]\leq \exp(-c n)$ (see, for example \cite{cannone_poisson}),
we infer that to show that $t(\perm{p}^{(n)})\geq \tfrac{1}{80}n^{2/3}$ a.a.s., it suffices to establish $t(\overline{S}_{n/2})\geq \tfrac{1}{80}n^{2/3}$ a.a.s.

Partition $[0,1]$ into $r\eqdef\lceil n^{2/3}\rceil$ equal intervals of length $1/r$ each, denoted $A_1,\dots,A_r$.
This induces a partition of $[0,1]^2$ into $r^2$ smaller squares of the form $A_i\times A_j$.
Sample a set $P$ from a Poisson process of intensity $n/2$ on $[0,1]^2$.
Make a bipartite graph $B$ whose parts are two copies of $[r]$, with $(i,j)$ being an edge if $A_i\times A_j$ contains at least two points of~$P$.
The edges are independent with probability $p=\Pr[\operatorname{Poisson}(n/2r^2)\geq 2]\geq \tfrac{1}{9}n^{-2/3}$, for large~$n$.
Clearly, every matching in $B$ corresponds to a pair of twins in the associated permutation.

\Cref{thm:random} follows once we show that $B$ is likely to contain a large matching. This
is well-known in the (very similar) context of the $G(n,p)$ model. We include such a proof for completeness.
\begin{claim}
Let $p\leq \tfrac{1}{6r}$. Then a random bipartite graph $G(r+r,p)$ contains a matching of size $pr^2/7$~a.a.s.
\end{claim}
\begin{proof}
Let $L\cup R$ be the bipartition. As long as $\abs{L}=\abs{R}\geq r/2$, do the following. 
Pick any vertex $v\in L$. It has a neighbor with probability $\geq p\abs{R}-p^2\binom{\abs{R}}{2}\geq pr/3$.
If $u\in R$ is a neighbor, match $u$ to $v$. Else, let $u$ be any vertex in $R$.
Remove $v$ from $L$ and $u$ from $R$. This way, we match $\operatorname{Binom}(r/2,pr/3)$
edges, which is at least $pr^2/7$ a.a.s.
\end{proof}

\textbf{Acknowledgment}. We thank Andrzej Dudek, Andrzej Ruci\'nski, Chengfei Xie, Zixiang Xu and two anonymous referees for feedback on earlier versions of this paper.\vspace*{-1.6ex}

\bibliographystyle{plain}
\bibliography{gawrontwins}

\end{document}